\documentclass[11pt]{article}

\usepackage[leqno]{amsmath}
\usepackage{amssymb}
\usepackage{amsthm}
\usepackage{amsfonts}
\usepackage{epsfig}
\usepackage{latexsym}
\usepackage{url}
\usepackage{hyperref}
\usepackage{ifthen}
\usepackage{comment}
\usepackage{subfigure}
\newcommand{\field}[1]{\mathbb{#1}}
\newcommand{\R}{\field{R}} 
\newcommand{\I}[1]{\mathbb{I}_{\left\{#1\right\}}} 
\newcommand{\tends}{{\rightarrow}} 
\newcommand{\ra}{{\rightarrow}} 
\newcommand{\subjectto}{\mbox{\rm subject to}} 
\renewcommand{\Re}{\R} 

\newcommand{\Cscr}{{\cal C}}

\newcommand{\Escr}{{\cal E}}

\newcommand{\Sscr}{{\cal S}}
\newcommand{\Tscr}{{\cal T}}

\newcommand{\Vscr}{{\cal V}}

\newcommand{\Xscr}{{\cal X}}

\newcommand\diag{\mathop{\mbox{{\rm diag}}}}

\newcommand\argmin{\mathop{\mbox{{\rm argmin}}}\limits}

\newcommand\minimize{\mathop{\mbox{{\rm minimize}}}\limits}
\newtheorem{theorem}{Theorem}

\newtheorem{lemma}{Lemma}
\newtheorem{assumption}{Assumption}
\newtheorem{definition}{Definition}

\numberwithin{equation}{section}
\newcommand{\emailhref}[1]{\href{mailto:#1}{\tt #1}} 
\title{Convergence of the \\ Min-Sum Algorithm for Convex Optimization}
\author{
Ciamac C. Moallemi \\
Electrical Engineering \\
Stanford University \\
email: \emailhref{ciamac@stanford.edu} \\
\and
Benjamin Van Roy \\
Management Science \&  Engineering \\
Electrical Engineering \\
Stanford University \\ 
email: \emailhref{bvr@stanford.edu}\\
}
\date{28 May 2007}

\begin{document}
\maketitle

\abstract{We establish that the min-sum message-passing algorithm and its asynchronous 
variants converge for a large class of unconstrained convex optimization problems.}

\section{Introduction}

Consider an optimization problem of the form
\begin{equation}\label{eq:separable}
\begin{array}{lll}
\minimize & F(x) = 
& \sum_{C \in \Cscr} f_C(x_C)
\\
\subjectto & & x\in\Xscr^V.
\end{array}
\end{equation}
Here, the vector of decision variables $x$ is indexed by a finite set
$V=\{1,\ldots,n\}$. Each decision variable takes values in the set
$\Xscr$. The set $\Cscr$ is a collection of subsets of the index set
$V$. This collection describes an additive decomposition of the
objective function. We associate with each set $C\in\Cscr$ a component
function (or {\it factor}) $f_C:\ \Xscr^C\ra\R$, which takes values as a
function of those components\footnote{Given a vector $x\in\Xscr^V$ and a
  subset $A\subset V$, we use the notation $x_A=(x_i,\ i\in A)\in\Xscr^A$
  for the vector of components of $x$ specified by the set $A$.} $x_C$
of the vector $x$ identified by the elements of $C$.

The min-sum algorithm is a method for optimization problems of the
form \eqref{eq:separable}. It is one of a class of methods know as
message-passing algorithms. These algorithms have been the subject of
considerable research recently across a number of fields, including
communications, artificial intelligence, statistical physics, and
theoretical computer science. Interest in message-passing algorithms
has been sparked by their success in solving certain classes of
NP-hard combinatorial optimization problems, such as the decoding of
low-density parity-check codes and turbo codes (e.g.,
\cite{Gallager63,Berrou93,Richardson01}), or the solution of certain
classes of satisfiability problems (e.g.,
\cite{Mezard02,Braunstein05}).

Despite their successes, message-passing algorithms remain poorly
understood.  For example, conditions for convergence and accurate
resulting solutions are not well characterized.

In this paper, we consider cases where $\Xscr=\R$, and the
optimization problem is continuous. One such case that has been
examined previously in the literature is where the objective is
pairwise separable (i.e., $|C|\leq 2$, for all $C\in\Cscr$) and the
component functions $\{ f_C(\cdot) \}$ are quadratic and convex. Here,
the min-sum algorithm is known to compute the optimal solution when it
converges \cite{Weiss01,Rusmevichientong01,Wainwright03}, and
sufficient conditions for convergence identify a broad class of
problems \cite{Moallemi06b,Malioutov06}.

Our main contribution is the analysis of cases where the functions are
convex but not necessarily quadratic. We establish that the min-sum
algorithm and its asynchronous variants converge for a large class of
such problems. The main sufficient condition is that of {\it scaled
  diagonal dominance}. This condition is similar to known sufficient
conditions for asynchronous convergence of other decentralized
optimization algorithms, such as coordinate descent and gradient
descent.

Analysis of the convex case has been an open challenge and its
resolution advances the state of understanding in the growing
literature on message-passing algorithms. Further, it builds a
bridge between this emerging research area and the better
established fields of convex analysis and optimization.

This paper is organized as follows.  The next section studies the
min-sum algorithm in the context of pairwise separable convex
programs, establishing convergence for a broad class of such problems.
Section~\ref{se:generalconvex} extends this result to more general
separable convex programs, where each factor can be a function of more
than two variables.  In Section~\ref{se:async}, we discuss how our
convergence results hold even with a totally asynchronous model of
computation.  When applied to a continuous optimization problem,
messages computed and stored by the min-sum algorithm are functions
over continuous domains.  Except in very special cases, this is not
feasible for digital computers, and in Section~\ref{se:impl}, we
discuss implementable approaches to approximating the behavior of the
min-sum algorithm.  We close by discussing possible extensions and
open issues in Section~\ref{se:open}.

\section{Pairwise Separable Convex Programs}

Consider first the case of pairwise separable programs. These
are programs of the form \eqref{eq:separable}, where $|C| \leq 2$, for
all $C\in \Cscr$. In this case, we can define an undirected graph
$(V,E)$ based on the objective function. This graph has a vertex set
$V$ corresponding to the decision variables, and an edge set $E$
defined by the pairwise factors,
\[
E = \{ C\in\Cscr:\ |C| = 2\}.
\]

\begin{definition}\label{def:pairwiseprog}
{\bf (Pairwise Separable Convex Program)} 
A pairwise separable convex program is an optimization problem of the form
\begin{equation}\label{eq:pairwiseopt}
\begin{array}{lll}
\minimize & F(x) = 
& \sum_{i\in V} f_i(x_i) + \sum_{(i,j)\in E} f_{ij}(x_i,x_j)
\\
\subjectto & & x\in\R^V,
\end{array}
\end{equation}
where the factors $\{ f_i(\cdot) \}$ are strictly
convex, coercive, and twice continuously differentiable, the
factors $\{f_{ij}(\cdot,\cdot) \}$ are convex and twice continuously
differentiable, and
\[
M \stackrel{\triangle}{=} \min_{i\in V}\ \inf_{x\in \R^V}\ \frac{\partial^2}{\partial x_i^2} F(x) > 0.
\]
\end{definition}

Under this definition, the objective function $F(x)$ is strictly
convex and coercive. Hence, we can define $x^*\in\R^V$ to be the
unique optimal solution.

\subsection{The Min-Sum Algorithm}

The min-sum algorithm attempts to minimize the objective function
$F(\cdot)$ by an iterative, message-passing procedure. For each vertex
$i\in V$, denote the set of neighbors of $i$ in the graph by
\[
N(i) = \{ j\in V:\ (i, j) \in E\}.
\]
Denote the set of edges with direction distinguished by
\[
\vec{E}=\{ (i,j)\in V\times V:\ i\in N(j) \}.
\]
At time $t$, each vertex $i$ keeps track of a ``message'' from each
neighbor $u\in N(i)$. This message takes the form of a function
$J^{(t)}_{u\rightarrow i}:\ \R\ra\R$. These incoming messages are
combined to compute new outgoing messages for each neighbor. The
message $J^{(t+1)}_{i\rightarrow j}(\cdot)$ from vertex $i$ to vertex
$j\in N(i)$ evolves according to
\begin{equation}\label{eq:Jupdate}
J_{i \rightarrow j}^{(t+1)}(x_j) = 
\min_{y_i} \left(f_i(y_i) + f_{ij}(y_i,x_j) 
+ \sum_{u \in N(i)\setminus j} J^{(t)}_{u \rightarrow i}(y_i)\right)
+ \kappa^{(t+1)}_{i\ra j}.
\end{equation}
Here, $\kappa^{(t+1)}_{i\ra j}$ represents an arbitrary offset term
that varies from message to message.  Only the relative values of the
function $J_{i \rightarrow j}^{(t+1)}(\cdot)$ matter, so the choice of
$\kappa^{(t+1)}_{i\ra j}$ does not influence relevant information.

At each time $t> 0$, a local objective function $b^{(t)}_i(\cdot)$ is
defined for each variable $x_i$ by
\begin{equation}
b^{(t)}_i(x_i) = f_i(x_i) + \sum_{u \in N(i)} J^{(t)}_{u\ra i}(x_i).
\end{equation}
An estimate $x^{(t)}_i$ can be obtained for the optimal value of the
variable $x_i$ by minimizing the local objective function:
\begin{equation}\label{eq:pairwiseest}
x^{(t)}_i = \argmin_{y_i} b^{(t)}_i(y_i).
\end{equation}

The min-sum algorithm requires an initial set of messages $\{
J^{(0)}_{i\ra j}(\cdot) \}$ at time $t=0$. We make the following
assumption regarding these messages:
\begin{assumption}\label{as:pairwiseinit}
  {\bf (Min-Sum Initialization)}
  Assume that the initial messages $\{ J^{(0)}_{i\ra j}(\cdot) \}$ are
  chosen to be twice continuously differentiable and so that, for each
  message $J^{(0)}_{i\ra j}(\cdot)$, there exists some $z_{i\ra
    j}\in\R$ with
\begin{equation}\label{eq:J0convex}
\frac{d^2}{d x_j^2} J^{(0)}_{i\ra j}(x_j) 
\geq
\frac{\partial^2}{\partial x_j^2} f_{ij}(z_{i\ra j}, x_j),
\quad\forall\ x_j\in\R.
\end{equation}
\end{assumption}

Assumption~\ref{as:pairwiseinit} guarantees that the messages at time
$t=0$ are convex functions. Examining the update equation
\eqref{eq:Jupdate}, it is clear that, by induction, this implies that
all future messages are also convex functions. Similarly, since the
functions $\{ f_i(\cdot) \}$ are strictly convex and coercive, and the
functions $\{ f_{ij}(\cdot,\cdot) \}$ are convex, it follows that the
optimization problem in the update equation \eqref{eq:Jupdate} is
well-defined and uniquely1 minimized. Finally, each local objective
function $b^{(t)}_i(\cdot)$ must strictly convex and coercive, and
hence each estimate $x^{(t)}_i$ is uniquely defined by
\eqref{eq:pairwiseest}.

Assumption~\ref{as:pairwiseinit} also requires that the initial
messages be sufficiently convex, in the sense of
\eqref{eq:J0convex}. As we will shortly demonstrate, this will be an
important condition for our convergence results. For the moment,
however, note that it is easy to select a set of initial messages
satisfying Assumption~\ref{as:pairwiseinit}. For example, one might
choose
\[
J^{(0)}_{i\ra j}(x_j) = f_{ij}(0,x_j).
\]

\subsection{Convergence}

Our goal is to understand conditions under which the min-sum algorithm
converges to the optimal solution $x^*$, i.e.
\[
\lim_{t\tends\infty} x^{(t)} = x^*.
\]

Consider the following diagonal dominance condition:
\begin{definition}{\bf (Scaled Diagonal Dominance)}
  An objective function $F:\ \R^V\ra\R$ is $(\lambda,w)$-scaled diagonally
  dominant if $\lambda$ is a scalar with $0 < \lambda < 1$ and
  $w\in\R^V$ is a vector with $w > 0$, so that for each $i\in V$ and
  all $x\in\R^V$,
\[
\sum_{j\in V\setminus i} 
w_j  \left|  \frac{\partial^2}{\partial x_i \partial x_j} F(x) \right|
\leq \lambda w_i \frac{\partial^2}{\partial x_i^2} F(x).
\]
\end{definition}

Our main convergence result is as follows:
\begin{theorem}\label{th:pairwiseconv}
  Consider a pairwise separable convex program with an objective
  function that is $(\lambda,w)$-scaled diagonally dominant. Assume
  that the min-sum algorithm is initialized in accordance with
  Assumption~\ref{as:pairwiseinit}.  Define the constant
\[
K = \frac{1}{M} \frac{\max_u w_u}{\min_u w_u}.
\]
Then, the iterates of the min-sum algorithm satisfy
\[
\| x^{(t)} - x^* \|_\infty
\leq
K
\frac{\lambda^t}{1 - \lambda}
\sum_{(u,v)\in\vec{E}}\ 
\left| 
\frac{d}{d x_v} J^{(0)}_{u\ra v}(x^*_v) 
-
\frac{\partial}{\partial x_v} f_{uv}(x^*_u, x^*_v)
\right|.
\]
Hence,
\[
\lim_{t\tends\infty}
x^{(t)}
=
x^*.
\]
\end{theorem}
\begin{proof}
The proof for Theorem~\ref{th:pairwiseconv} will be provided in
Section~\ref{se:convproof}. 
\end{proof}

We can compare Theorem~\ref{th:pairwiseconv} to existing results on
min-sum convergence in the case of where the objective function
$F(\cdot)$ is quadratic. Rusmevichientong and Van Roy
\cite{Rusmevichientong01} developed abstract conditions for
convergence, but these conditions are difficult to verify in practical
instances. Convergence has also been established in special cases
arising in certain applications \cite{Moallemi06a,Montanari05}.

More closely related to our current work, Weiss and Freeman
\cite{Weiss01} established convergence when the factors $\{
f_i(\cdot), f_{ij}(\cdot,\cdot) \}$ are quadratic, the single-variable
factors $\{ f_i(\cdot) \}$ are strictly convex, and the pairwise
factors $\{ f_{ij}(\cdot,\cdot) \}$ are convex and diagonally
dominated, i.e.
\[
\left|  \frac{\partial^2}{\partial x_i \partial x_j} f_{ij}(x_i,x_j) \right|
\leq  \frac{\partial^2}{\partial x_i^2} f_{ij}(x_i,x_j),
\quad\forall\ (i,j)\in E,\ x_j,x_j\in\R.
\]
The results of Malioutov, et al. \cite{Malioutov06} and our prior work
\cite{Moallemi06b} remove the diagonal dominance assumption. However,
all of these results are special cases of
Theorem~\ref{th:pairwiseconv}. In particular, if the a quadratic
objective function $F(\cdot)$ decomposes into pairwise factors so that the
single-variable factors are quadratic and strictly convex, and the
pairwise factors are quadratic convex, then $F(\cdot)$ must be scaled
diagonally dominant. This can be established as a consequence of the
Perron-Frobenius theorem \cite{Malioutov06}. Finally, as we will see
in Section~\ref{se:generalconvex}, Theorem~\ref{th:pairwiseconv} also
generalizes beyond pairwise decompositions.

\subsection{The Computation Tree}\label{se:comptree}

In order to prove Theorem~\ref{th:pairwiseconv}, we first introduce
the notion of the {\it computation tree}. This is a useful device in
the analysis of message-passing algorithms, originally introduced by
Wiberg \cite{Wiberg96}.  Given a vertex $r\in V$ and a time $t$, the
computation tree defines an optimization problem that is constructed
by ``unrolling'' all the optimizations involved in the computation of
the min-sum estimate $x^{(t)}_r$.

Formally, the computation tree is a graph $\Tscr=(\Vscr,\Escr)$ where
each vertex $i\in\Vscr$ is in labeled by a vertex $\tilde{i}\in V$ in
the original graph, through a mapping $\sigma:\ \Vscr\ra V$. This
mapping is required to preserve the edge structure of the graph, so
that if $(i,j)\in\Escr$, then $(\sigma_i, \sigma_j)\in E$. Given a
vertex $i\in \Vscr$, we will abuse notation and refer to the
corresponding vertex $\sigma_i\in V$ in the original graph simply by
$i$.

\begin{figure}[htpb]
\centering
\subfigure{\includegraphics{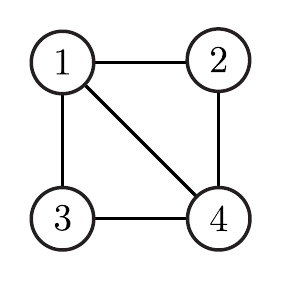}}
\subfigure{\includegraphics{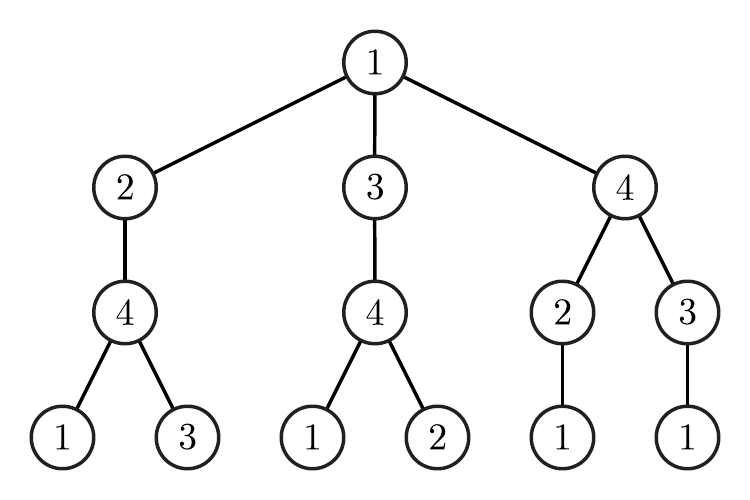}}
\caption{A graph and the corresponding computation tree, rooted at
  vertex $1$ and of depth $t=3$. The vertices in the computation tree are
  labeled according to the corresponding vertex in the original
  graph.}
\label{fig:pairwise-tree}
\end{figure}

Fixing a vertex $r\in V$ and a time $t$, the computation tree rooted
at $r$ and of depth $t$ is defined in an iterative fashion. Initially,
the tree consists of a root single vertex corresponding to $r$. At
each subsequent step, the leaves in the computation tree are
examined. Given a leaf $i$ with a parent $j$, a vertex $u$ and an edge
$(u,i)$ are added to the computation tree corresponding to each
neighbor of $i$ excluding $j$ in the original graph. This process is
repeated for $t$ steps. An example of the resulting graph is
illustrated in Figure~\ref{fig:pairwise-tree}.

Given the graph $\Tscr=(\Vscr,\Escr)$, and the correspondence mapping
$\sigma$, define a decision variable $x_i$ for each vertex $i\in
\Vscr$. Define a pairwise separable objective function $F_\Tscr:\
\R^\Vscr\ra\R$, by considering factors of the form:
\begin{enumerate}
\item For each $i\in\Vscr$, add a single-variable factor $f_i(x_i)$ by
  setting $f_i(x_i) \stackrel{\triangle}{=} f_{\sigma_i}(x_i)$.
\item For each $(i,j)\in\Vscr$, add a pairwise factor
  $f_{ij}(x_i,x_j)$ by setting $f_{ij}(x_i,x_j)
  \stackrel{\triangle}{=} f_{\sigma_i \sigma_j}(x_i,x_j)$.
\item For each $i\in \Vscr$ that is a leaf vertex with parent $j$, add a
  single-variable factor $J^{(0)}_{u\ra \sigma_i}(x_i)$, for each
  neighbor $u\in N(\sigma_i)\setminus\sigma_j$ of $i$ in the original
  graph, excluding $j$.
\end{enumerate}

Now, let $\tilde{x}$ be the optimal solution to the minimization of
the computation tree objective $F_\Tscr(\cdot)$. By inductively
examining the operation of the min-sum algorithm, it is easy to
establish that the component $\tilde{x}_r$ of this solution at the
root of the tree is precisely the min-sum estimate $x^{(t)}_r$.

The following lemma establishes that the computation tree inherits the
scaled diagonal dominance property from the original objective
function.
\begin{lemma}\label{le:comptreedd}
  Consider a pairwise separable convex program with an objective
  function that is $(\lambda,w)$-scaled diagonally dominant.  Assume
  that the min-sum algorithm is initialized in accordance with
  Assumption~\ref{as:pairwiseinit}, and let $\Tscr=(\Vscr,\Escr)$ be a
  computation tree associated with this program. Then, the computation
  tree objective function $F_\Tscr(\cdot)$ is also
  ($\lambda,w)$-scaled diagonally dominant.
\end{lemma}
\begin{proof}
  Given a vertex $i\in \Vscr$, let $N^\Vscr(i)$ be the neighborhood in
  the computation tree, and let $N(i)$ be the neighborhood of the
  corresponding vertex in the original graph.  If $i\in \Vscr$ is an
  interior vertex of the computation tree, then
\[
\begin{split}
\lefteqn{
\sum_{u\in \Vscr\setminus i} 
w_u  \left|  \frac{\partial^2}{\partial x_i \partial x_u} F_\Tscr(x) \right|
} & \\
& =
\sum_{u\in N^\Vscr(i)}
w_u  \left|  \frac{\partial^2}{\partial x_i \partial x_u} f_{iu}(x_i,x_u) \right|
\\
& \leq
 \lambda w_i\left(
\frac{\partial^2}{\partial x_i^2} f_i(x_i)
+
 \sum_{u\in N^\Vscr(i)} \frac{\partial^2}{\partial x_i^2} f_{iu}(x_i,x_u)\right)
\\
& =
\lambda w_i \frac{\partial^2}{\partial x_i^2} F_\Tscr(x),
\end{split}
\]
where the inequality follows from the scaled diagonal dominance of the
original objective function $F(\cdot)$.

Similarly, if $i$ is a leaf vertex with parent $j$,
\[
\begin{split}
\lefteqn{
\sum_{u\in \Vscr\setminus i} 
w_u  \left|  \frac{\partial^2}{\partial x_i \partial x_u} F_\Tscr(x) \right|
} & \\
& =
w_j  \left|  \frac{\partial^2}{\partial x_i \partial x_j} f_{ij}(x_i,x_j) \right|
\\
&
\leq
w_j  \left|  \frac{\partial^2}{\partial x_i \partial x_j} f_{ij}(x_i,x_j) \right|
+ \sum_{u\in N(i)\setminus j} 
w_u  \left|  \frac{\partial^2}{\partial x_i \partial x_u} f_{iu}(x_i,z_{u\ra i}) \right|
\\
&
\leq
\lambda w_i 
\left(
\frac{\partial^2}{\partial x_i^2} f_i(x_i)
+
\frac{\partial^2}{\partial x_i^2} f_{ij}(x_i,x_j)
+
 \sum_{u\in N(i)\setminus j}
\frac{\partial^2}{\partial x_i^2} f_{iu}(x_i,z_{u\ra i})
\right)
\\
&
\leq
\lambda w_i 
\left(
\frac{\partial^2}{\partial x_i^2} f_i(x_i)
+
\frac{\partial^2}{\partial x_i^2} f_{ij}(x_i,x_j)
+
 \sum_{u\in N(i)\setminus j}
 \frac{\partial^2}{\partial x_i^2}
 J^{(0)}_{u\ra i}(x_i)
\right)
\\
& =
\lambda w_i \frac{\partial^2}{\partial x_i^2} F_\Tscr(x).
\end{split}
\]
Here, the second inequality follows from the scaled diagonal dominance
of the original objective function $F(\cdot)$, and the third
inequality follows from Assumption~\ref{as:pairwiseinit}.
\end{proof}

\subsection{Proof of Theorem~\ref{th:pairwiseconv}}\label{se:convproof}

In order to prove Theorem~\ref{th:pairwiseconv}, we will study the
evolution of the min-sum algorithm under a set of linear
perturbations. Consider an arbitrary vector $p\in\R^{\vec{E}}$ with
one component $p_{i\ra j}$ for each $i\in V$ and $j\in N(i)$.  Given
an arbitrary vector $p$, define $\{ J^{(t)}_{i\ra j}(\cdot,p) \}$ to
be the set of messages that evolve according to
\begin{equation}\label{eq:Jpstarupdate}
\begin{split}
J^{(0)}_{i\ra j}(x_j, p) & = 
J^{(0)}_{i\ra j}(x_j)
+ p_{i\ra j}  x_j,
\\
J_{i \rightarrow j}^{(t+1)}(x_j,p) & = 
\min_{y_i} \left(f_i(y_i) + f_{ij}(y_i,x_j) 
+ \sum_{u \in N(i)\setminus j} J^{(t)}_{u \rightarrow i}(y_i,p)\right)
\\
& \quad
+ \kappa^{(t+1)}_{i\ra j}.
\end{split}
\end{equation}
Similarly, define $\{ b^{(t)}_i(\cdot,p) \}$ and $\{x^{(t)}_i(p)\}$ to be
the resulting local objective functions and optimal value estimates
under this perturbation:
\begin{align*}
b^{(t)}_i(x_i,p) & = f_i(x_i) + \sum_{u \in N(i)} J^{(t)}_{u\ra i}(x_i,p),
\\
x^{(t)}_i(p) & = \argmin_{y_i} b^{(t)}_i(y_i,p).
\end{align*}

The following simple lemma gives a particular choice of $p$ for which
the min-sun algorithm yields the optimal solution at every time.
\begin{lemma}\label{le:pstarexact}
  Define the vector $p^*\in\R^{\vec{E}}$ by setting, for each $i \in
  V$ and $j\in N(i)$,
\[
p^*_{i\ra j} =
\frac{\partial}{\partial x_j} f_{ij}(x^*_i, x^*_j)
-
\frac{d}{d x_j} J^{(0)}_{i\ra j}(x^*_j).
\]
Then, at every time $t \geq 0$,
\begin{equation}\label{eq:Jpstar}
\frac{\partial}{\partial x_j}
J^{(t)}_{i\ra j}(x^*_j, p^*) = 
\frac{\partial}{\partial x_j} f_{ij}(x^*_i, x^*_j),
\end{equation}
and $x^{(t)}_j(p^*) = x^*_j$.
\end{lemma}
\begin{proof}
Note that the first order optimality conditions for $F(x)$ at $x^*$ imply that, for each $j\in V$, 
\[
\frac{d}{d x_j} f_i(x^*_j) + \sum_{i\in N(j)} 
\frac{\partial}{\partial x_j} f_{ij}(x^*_i, x^*_j) = 0.
\]
If \eqref{eq:Jpstar} holds at time $t$, this is exactly the first
order optimality condition for the minimization of
$b^{(t)}_j(\cdot,p^*)$, thus $x^{(t)}_j(p^*) = x^*_j$.

Clearly \eqref{eq:Jpstar} holds at time $t=0$. Assume it holds at time
$t \geq 0$. Then, when $x_j=x^*_j$, the minimizing value of $y_i$ in \eqref{eq:Jpstar} is $x^*_i$. Hence, \eqref{eq:Jpstar} holds at time $t+1$.
\end{proof}

Next, we will bound the sensitivity of the estimate $x^{(t)}_i(p)$ to
the choice of $p$. The main technique employed here is analysis of the
computation tree described in Section~\ref{se:comptree}. In
particular, the perturbation $p$ impacts the computation tree only
through the leaf vertices at depth $t$. The scaled diagonal dominance
property of the computation tree, provided by
Lemma~\ref{le:comptreedd}, can then be used to guarantee that this
impact is diminishing in $t$.

\begin{lemma}\label{le:psensitivity}
  We have, for all $p\in\R^{\vec{E}}$, $r \in V$, $(u,v)\in\vec{E}$,
  and $t \geq 0$,
\[
\left|
\frac{\partial}{\partial p_{u\ra v}} x^{(t)}_r(p)
\right|
\leq K \frac{\lambda^t}{1-\lambda}.
\]
\end{lemma}
\begin{proof}
  Fix $r \in V$, and let $\Tscr=(\Vscr,\Escr)$ be the computation tree
  rooted at $r$ after $t$ time steps. Let $F_\Tscr(x,p)$ be the
  objective value of this computation tree, and let
\[
\tilde{x}(p) = \argmin_x F_\Tscr(x,p),
\]
so that
\[
\tilde{x}_r(p) = x^{(t)}_r(p).
\]
By the first order optimality conditions, for any $j\in \Vscr$,
\[
\frac{\partial}{\partial x_j} F_\Tscr(\tilde{x}(p),p) = 0.
\]
If $j$ is an interior vertex of $\Tscr$, this becomes
\begin{equation}\label{eq:fo1}
\frac{d}{d x_j} f_j(\tilde{x}_j(p))
+ 
\sum_{i \in N(j)} 
\frac{\partial}{\partial x_j} f_{ij}(\tilde{x}_i(p),\tilde{x}_j(p)) = 0.
\end{equation}
If $j$ is a leaf with parent $u$, we have
\begin{multline}\label{eq:fo2}
\frac{d}{d x_j} f_j(\tilde{x}_j(p))
+
\frac{\partial}{\partial x_j} f_{uj}(\tilde{x}_u(p),\tilde{x}_j(p))
\\
+ 
\sum_{i \in N(j)\setminus u} 
\left(\frac{\partial}{\partial x_j} J^{(0)}_{i\ra j}(\tilde{x}_j(p))  
+ p_{i\ra j}\right) = 0.
\end{multline}
Now, fixed some directed edge $(a,b)$, and differentiate
\eqref{eq:fo1}--\eqref{eq:fo2} with respect to $p_{a\ra b}$. We have,
for an interior vertex $j$,
\[
\begin{split}
0 & = 
\frac{d^2}{d x_j^2} f_j(\tilde{x}_j(p))
\frac{\partial}{\partial p_{a\ra b}} \tilde{x}_j(p)
\\
& \quad
+
\sum_{i \in N(j)} 
\frac{\partial^2}{\partial x_j^2} f_{ij}(\tilde{x}_i(p),\tilde{x}_j(p)) 
\frac{\partial}{\partial p_{a\ra b}} \tilde{x}_j(p)
\\
& \quad
+
\sum_{i \in N(j)} 
\frac{\partial^2}{\partial x_i \partial x_j} f_{ij}(\tilde{x}_i(p),\tilde{x}_j(p)) 
\frac{\partial}{\partial p_{a\ra b}} \tilde{x}_i(p),
\end{split}
\]
and for a leaf vertex $j$ with parent $u$,
\[
\begin{split}
0 & =
\frac{d^2}{d x_j^2} f_j(\tilde{x}_j(p))
\frac{\partial}{\partial p_{a\ra b}} \tilde{x}_j(p)
\\
& \quad
+
\frac{\partial^2}{\partial x_j^2} f_{uj}(\tilde{x}_u(p),\tilde{x}_j(p))
\frac{\partial}{\partial p_{a\ra b}} \tilde{x}_j(p)
\\
& \quad
+
\frac{\partial^2}{\partial x_u \partial x_j} f_{uj}(\tilde{x}_u(p),\tilde{x}_j(p))
\frac{\partial}{\partial p_{a\ra b}} \tilde{x}_u(p)
\\
& \quad
+ 
\sum_{i \in N(j)\setminus u} 
\left(
\frac{\partial^2}{\partial x_j^2} J^{(0)}_{i\ra j}(\tilde{x}_j(p)) 
\frac{\partial}{\partial p_{a\ra b}} \tilde{x}_j(p)
+
\I{(a,b) = (i,j)}\right).
\end{split}
\]
We can write this system of equations in matrix form, as
\begin{equation}\label{eq:lineq}
\Gamma v^{a\ra b}  + h^{a\ra b} = 0.
\end{equation}
Here, $v^{a\ra b}\in\R^\Vscr$ is a vector with components
\[
v^{a\ra b}_j = \frac{\partial}{\partial p_{a\ra b}} \tilde{x}_j(p).
\]
The vector $h^{a\ra b}\in\R^\Vscr$ has components
\[
h^{a\ra b}_j = 
\I{\text{$j$ is a leaf vertex of type $a$ with a parent of type $b$}}.
\]
The symmetric matrix $\Gamma \in \R^{\Vscr\times\Vscr}$ has components
as follows:
\begin{enumerate}
\item If $j$ is an interior vertex,
\[
\Gamma_{j j} =
\frac{d^2}{d x_j^2} f_j(\tilde{x}_j(p))
+
\sum_{i \in N(j)}
\frac{\partial^2}{\partial x_j^2} f_{ij}(\tilde{x}_i(p),\tilde{x}_j(p)).
\]
\item If $j$ is an interior vertex and $i\in N(j)$,
\[
\Gamma_{i j} =
\frac{\partial^2}{\partial x_i \partial x_j} f_{i j}(\tilde{x}_i(p),\tilde{x}_j(p)).
\]
\item If $j$ is a leaf vertex with parent $u$,
\[
\begin{split}
\Gamma_{jj} & = 
\frac{d^2}{d x_j^2} f_j(\tilde{x}_j(p))
+
\frac{\partial^2}{\partial x_j^2} f_{uj}(\tilde{x}_u(p),\tilde{x}_j(p))
\\
& \quad +
\sum_{i \in N(j)\setminus u} 
\frac{\partial^2}{\partial x_j^2} J^{(0)}_{i\ra j}(\tilde{x}_j(p)),
\\
\Gamma_{uj} & = 
\frac{\partial^2}{\partial x_u \partial x_j} 
f_{uj}(\tilde{x}_u(p),\tilde{x}_j(p)).
\end{split}
\]
\item All other entries of $\Gamma$ are zero.
\end{enumerate}

Note that $\Gamma=\nabla^2_x F_\Tscr(\tilde{x}(p), p)$. Then,
Lemma~\ref{le:comptreedd} implies that
\begin{equation}\label{eq:gammacont}
\sum_{i\in\Vscr\setminus j} w_i |\Gamma_{ij}|
\leq \lambda w_j \Gamma_{jj}.
\end{equation}

Define, for vectors $x \in \R^\Vscr$, the weighted sup-norm
\[
\| x \|^w_\infty = \max_{j\in \Vscr} |x_j|/w_j.
\]
For a linear operator $A:\ \R^\Vscr\ra\R^\Vscr$, The corresponding
induced operator norm is given by
\[
\| A \|^w_\infty
= \max_{j \in \Vscr} \frac{1}{w_j} \sum_{i\in \Vscr} w_i |A_{ji}|.
\]

Define the matrices
\[
\begin{split}
D & = \diag( \Gamma ), \\
R & = I - D^{-1} \Gamma.
\end{split}
\]
Then, \eqref{eq:gammacont} implies that
\[
\| R \|^w_\infty \leq \lambda < 1.
\]
Hence, the matrix $I - R = D^{-1} \Gamma$ is invertible, and
\[
\left(D^{-1} \Gamma\right)^{-1} = (I - R)^{-1} = \sum_{s=0}^\infty R^s.
\]

Examining the linear equation \eqref{eq:lineq}, we have
\[
v^{a\ra b}
= - \Gamma^{-1} h^{a\ra b}
= - (I - R)^{-1} D^{-1} h^{a\ra b}
= - \sum_{s=0}^\infty R^s D^{-1} h^{a \ra b}.
\]
We are interested in bounding the value of the component $v^{a\ra
  b}_r$ (recall that $v^{a\ra b}_r = \partial x^{(t)}_r(p)/\partial
p_{a \ra b}$). Hence, we have
\[
v^{a\ra b}_r = - \sum_{s=0}^\infty \left[R^s D^{-1} h^{a \ra b}\right]_r.
\]
Since $h^{a\ra b}$ is zero on interior vertices, and any leaf vertex is
distance $t$ from the root $r$, we have
\[
\left[R^s D^{-1} h^{a \ra b}\right]_r = 0,\quad\forall\ s < t.
\]
Thus,
\[
v^{a\ra b}_r = - \sum_{s=t}^\infty \left[R^s D^{-1} h^{a \ra b}\right]_r.
\]
Then,
\[
\begin{split}
|v^{a\ra b}_r| / w_r  
& \leq 
\left\|
\sum_{s=t}^\infty R^s D^{-1} h^{a \ra b} \right\|^w_\infty
\\
& \leq 
\sum_{s=t}^\infty 
\left\| R^s \right\|^w_\infty \left\| D^{-1} h^{a \ra b} \right\|^w_\infty
\\
& \leq
\frac{\lambda^t}{1 - \lambda} \left\| D^{-1} h^{a \ra b} \right\|^w_\infty
\\
& \leq
\frac{\lambda^t}{1 - \lambda}
\max_{i \in V} \sup_{x\in\R^V} \left( w_i \frac{\partial^2}{\partial x_i^2} F(x) \right)^{-1}
\\
& \leq
M \frac{\lambda^t}{1 - \lambda} \max_{i\in \Vscr} \frac{1}{w_i}.
\end{split}
\]
\end{proof}

The following lemma combines the results from
Lemmas~\ref{le:pstarexact} and \ref{le:psensitivity}.
Theorem~\ref{th:pairwiseconv} follows by taking $p=0$.

\begin{lemma}\label{le:xconv}
Given an arbitrary vector $p\in\R^{\vec{E}}$,
\[
\| x^{(t)}(p) - x^* \|_\infty
\leq
K
\frac{\lambda^t}{1 - \lambda}
\sum_{(u,v)\in\vec{E}}\ 
\left| 
p_{u\ra v}
-
p^*_{u\ra v}
\right|.
\]
\end{lemma}
\begin{proof}
For any $j\in V$, define
\[
g_j^{(t)}(\theta) = x^{(t)}_j(\theta p + (1 - \theta) p^*).
\]
We have, from Lemma~\ref{le:pstarexact},
\[
x^{(t)}_j(p) - x^*_j = x^{(t)}_j(p) - x^{(t)}_j(p^*) = g^{(t)}_j(1) - g^{(t)}_j(0).
\]
By the mean value theorem and Lemma~\ref{le:psensitivity},
\[
\begin{split}
\lefteqn{
|
x^{(t)}_j(p) - x^*_j
|
} & \\
& 
\leq
\sup_{\theta \in [0,1]}
\left| \frac{d}{d\theta} g^{(t)}_j(\theta) \right|
\\
& \leq
\sup_{\theta \in [0,1]}
\sum_{(u,v)\in\vec{E}}
\left| 
\frac{\partial}{\partial p_{u\ra v}}
x^{(t)}_j(\theta p + (1 - \theta) p^*)
\right|
| p_{u \ra v} - p^*_{u \ra v} |
\\
& \leq
K \frac{\lambda^t}{1-\lambda}
\sum_{(u,v)\in\vec{E}}
| p_{u \ra v} - p^*_{u \ra v} |.
\end{split}
\]
\end{proof}

\section{General Separable Convex Programs}\label{se:generalconvex}

In this section we will  consider convergence of the min-sum algorithm
for more general separable  convex programs. In particular, consider a
vector of  real-valued decision variables  $x \in \R^V$, indexed  by a
finite set $V$, and a hypergraph $(V,\Cscr)$, where the set $\Cscr$ is
a collection of subsets (or ``hyperedges'') of the vertex set $V$.

\begin{definition}{\bf (General Separable Convex Program)}
A general separable convex program 
is an optimization problem of the form
\begin{equation}\label{eq:generalopt}
\begin{array}{lll}
\minimize & F(x) = 
& \sum_{i\in V} f_i(x_i) + \sum_{C\in \Cscr} f_{C}(x_C)
\\
\subjectto & & x\in\R^V,
\end{array}
\end{equation}
where the factors $\{ f_i(\cdot) \}$ are strictly
convex, coercive, and twice continuously differentiable, the
factors $\{f_{C}(\cdot) \}$ are convex and twice continuously
differentiable, and
\[
M \stackrel{\triangle}{=} \min_{i\in V}\ \inf_{x\in \R^V}\ \frac{\partial^2}{\partial x_i^2} F(x) > 0.
\]
\end{definition}

In this setting, the min-sum algorithm operates by passing messages
between vertices and hyperedges. In particular, denote the set of neighbor hyperedges to a vertex $i\in V$ by
\[
N_f(i) = \{ C\in\Cscr\ :\ i\in C\},
\]
The min-sum update equations take the form
\begin{equation}\label{eq:generalJupdate}
\begin{split}
J^{(t+1)}_{i\ra C}(x_i)
&
=
f_i(x_i) + 
\sum_{C'\in N_f(i) \setminus C} J^{(t)}_{C'\ra i}(x_i) + \kappa^{(t+1)}_{i \ra C},
\\
J^{(t+1)}_{C \ra i}(x_i) 
& 
= 
\minimize_{y_{C\setminus i}}\  f_C(x_i,y_{C\setminus i}) + \sum_{i'\in C\setminus i} 
J^{(t+1)}_{i'\ra C}(y_{i'}) + \kappa^{(t+1)}_{C\ra i}.
\end{split}
\end{equation}
Local objective functions and estimates of the optimal solution are defined by
\begin{align*}
b^{(t)}_i(x_i) & = f_i(x_i) + \sum_{C \in N_f(i)} J^{(t)}_{C\ra i}(x_i), \\
x^{(t)}_i & = \argmin_{y_i} b^{(t)}_i(y_i).
\end{align*}

We will make the following assumption on the initial messages:
\begin{assumption}\label{as:generalinit} {\bf (Min-Sum Initialization)}
  Assume that the initial messages $\{ J^{(0)}_{C\ra j}(\cdot) \}$ are
  chosen to be twice continuously differentiable and so that, for each
  message $J^{(0)}_{C\ra j}(\cdot)$, there exists some $z_{C\ra
    j}\in\R^{C\setminus i}$ with
\[
\frac{d^2}{d x_j^2} J^{(0)}_{C\ra j}(x_j) 
\geq
\frac{\partial^2}{\partial x_j^2} f_{C}(x_j,z_{C\ra j}),
\quad\forall\ x_j\in\R.
\]
\end{assumption}

Then, we have the following analog of Theorem~\ref{th:pairwiseconv}:
\begin{theorem}\label{th:generalconv}
  Consider a general separable convex program. Assume that either:
\begin{itemize}
\item[(i)] The objective function $F(x)$ is scaled diagonally
  dominant, and each pair of vertices $i, j\in V$ participate in at
  most one common factor. That is,
\[
| \{ C \in \Cscr\ :\ (i,j)\subset C\} | \leq 1,\quad\forall\ i,j\in V.
\]
\item[(ii)] The factors $\{ f_C(\cdot) \}$ are individually scaled
  diagonally dominant, in the sense that exists a scalar $\lambda\in
  (0,1)$ and a vector $w\in\R^V$, with $w > 0$, so that for all $C \in
  \Cscr$, $i\in C$, and $x_C\in\R^C$,
\[
\sum_{j\in C\setminus i} 
w_j  \left|  \frac{\partial^2}{\partial x_i \partial x_j} f_C(x_C) \right|
\leq \lambda w_i \frac{\partial^2}{\partial x_i^2} f_C(x_C).
\]
\end{itemize}
Assume that the min-sum algorithm is initialized in accordance with
Assumption~\ref{as:generalinit}. Define the constant
\[
K = \frac{1}{M} \frac{\max_u w_u}{\min_u w_u}.
\]
Then, the iterates of the min-sum algorithm satisfy
\[
\| x^{(t)} - x^* \|_\infty
\leq
K
\frac{\lambda^t}{1 - \lambda}
\sum_{C\in\Cscr}
\sum_{v \in C}\ 
\left| 
\frac{d}{d x_v} J^{(0)}_{C\ra v}(x^*_v) 
-
\frac{\partial}{\partial x_v} f_{C}(x^*_C)
\right|.
\]
Hence,
\[
\lim_{t\tends\infty}
x^{(t)}
=
x^*.
\]
\end{theorem}
\begin{proof}
  This result can be proved using the same method as
  Theorem~\ref{th:pairwiseconv}. The main modification required is the
  development of a suitable analog of Lemma~\ref{le:comptreedd}. In
  the general case, scaled diagonal dominance of the computation tree
  does not follow from scaled diagonal dominance of the objective
  function $F(x)$. However, it is easy to verify that either of the
  hypotheses (i) or (ii) imply scaled diagonal dominance of the
  computation tree. The balance of the proof proceeds as in
  Section~\ref{se:convproof}.
\end{proof}

\section{Asynchronous Convergence}\label{se:async}

The convergence results of Theorems~\ref{th:pairwiseconv} and
\ref{th:generalconv} assumed a synchronous model of computation. That
is, each message is updated at every time step in parallel. The
min-sum update equations \eqref{eq:Jupdate} and
\eqref{eq:generalJupdate} are naturally decentralized, however. If we
consider the application of the min-sum algorithm in distributed
contexts, it is necessary to consider convergence under an {\it
  asynchronous} model of computation. In this section, we will
establish that Theorems~~\ref{th:pairwiseconv} and
\ref{th:generalconv} extend to an asynchronous setting.

Without loss of generality, consider the pairwise case.  Assume that
there is a processor associated with each vertex $i$ in the graph, and
that this processor is responsible for computing the message $J_{i\ra
  j}(\cdot)$, for each neighbor $j$ of vertex $i$. Each processor
occasionally communicates its messages to neighboring processors, and
occasionally computes new messages based on the most recent messages
it has received. Define the $T^{i}$ to be the set of times at which
new messages are computed. Define $0 \leq \tau_{j\ra i}(t) \leq t$ to
be the last time the processor at vertex $j$ communicated to the
processor at vertex $i$. Then, the messages evolve according to
\[
\begin{split}
J_{i \rightarrow j}^{(t+1)}(x_j) & = 
\min_{y_i} \left(f_i(y_i) + f_{ij}(y_i,x_j) 
+ \sum_{u \in N(i)\setminus j} J^{(\tau_{u\ra i}(t))}_{u \rightarrow i}(y_i)\right)
\\
&\quad\quad
+ \kappa^{(t+1)}_{i\ra j},
\end{split}
\]
if $t\in T^i$, and 
\[
J_{i \rightarrow j}^{(t+1)}(x_j) = 
J_{i \rightarrow j}^{(t)}(x_j),
\]
otherwise.

We will make the following assumption \cite{BertsekasPDP}:
\begin{assumption} {\bf (Total Asynchronism)}
Assume that 
each set $T^i$ is infinite, and
that if $\{t_k\}$ is a sequence in $T^i$ tending to infinity, then
\[
\lim_{k\rightarrow \infty} \tau_{i\rightarrow j}(t_k) = \infty,
\]
for each neighbor $j \in N(i)$.
\end{assumption}
Total asynchronism is a very mild assumption. It guarantees that each
component is updated infinitely often, and that processors eventually
communicate with neighboring processors. It allows for arbitrary
delays in communication, and even the out-of-order arrival of messages
between processors.

Theorem~\ref{th:pairwiseconv} can be extended to the totally
asynchronous setting. To see this, note that we can repeat the
construction of the computation tree in Section~\ref{se:comptree}. As
in the synchronous case, the initial messages only impact the leaves
of computation tree. The total asynchronism assumption guarantees that
these leaves are, eventually, arbitrarily far away from the root of
the computation tree. The arguments in Lemma~\ref{le:psensitivity}
then imply that the optimal value at the root of the computation
tree is insensitive to the choice of initial messages. Convergence
follows, as in Section~\ref{se:convproof}.

The scaled diagonal dominance requirement of our convergence result is
similar to conditions required for the totally asynchronous
convergence of other optimization algorithms. Consider, for example, a
decentralized coordinate descent algorithm. Here, the processor
associated with vertex $i$ maintains an estimate $x^{(t)}_i$ of the
$i$th component of the optimal solution at time $t$. These estimates
are updated according to
\[
x^{(t+1)}_i =
\argmin_{y_i} f_i(y_i) + 
\sum_{u\in N(i)} f_{oi}(x^{(\tau_{u\ra i}(t))}_u,y_i),
\]
if $t \in T^i$, and $x^{(t+1)}_i = x^{(t)}_i$, otherwise.

Similarly, consider a decentralized gradient method, where
\[
x^{(t+1)}_i =
x^{(t)}_i - \alpha 
\frac{\partial}{\partial x_i} 
\left(
 f_i\left(x^{(t)}_i\right) + 
\sum_{u\in N(i)}
f_{ui}\left(x^{(\tau_{u\ra i}(t))}_u,x^{(t)}_i\right)
\right),
\]
if $t\in T^i$, and $x^{(t+1)}_i = x^{(t)}_i$, otherwise, for some
small positive step size $\alpha$.  These methods are not guaranteed
to converge for arbitrary pairwise separable convex optimization
problems. Typically, some sort of diagonal dominance condition is
needed \cite{BertsekasPDP}.

\section{Implementation}\label{se:impl}

The convergence theory we have presented elucidates properties of the
min-sum algorithm and builds a bridge to the more established areas of
convex analysis and optimization.  However, except in very special
cases, the algorithm as we have formulated it can not be implemented
on a digital computer because the messages that are computed and
stored are functions over continuous domains.  In this section, we
present two variations that can be implemented to approximate behavior
of the min-sum algorithm.  For simplicity, we restrict attention to
the case of synchronous min-sum for pairwise separable convex
programs.

Our first approach approximates messages using quadratic functions and
can be viewed as a hybrid between the min-sum algorithm and Newton's
method.  It is easy to show that, if the single-variable factors
$\{f_i(\cdot)\}$ are positive definite quadratics and the pairwise
factors $\{f_{ij}(\cdot,\cdot)\}$ are positive semidefinite
quadratics, then min-sum updates map quadratic messages to quadratic
messages.  The algorithm we propose here maintains a running estimate
$\tilde{x}^{(t)}$ of the optimal solution, and at each time
approximates each factor by a second-order Taylor expansion.  In
particular, let $\tilde{f}_i^{(t)}(\cdot)$ be the second-order Taylor
expansion of $f_i(\cdot)$ around $\tilde{x}_i^{(t)}$ and let
$\tilde{f}_{ij}^{(t)}(\cdot,\cdot)$ be the second-order Taylor
expansion of $f_{ij}(\cdot,\cdot)$ around $(\tilde{x}_i^{(t)},
\tilde{x}_j^{(j)})$.  Quadratic messages are updated according to
\begin{equation}\label{eq:quadratic-update}
J_{i \rightarrow j}^{(t+1)}(x_j) = 
\min_{y_i} \left(\tilde{f}^{(t)}_i(y_i) + \tilde{f}^{(t)}_{ij}(y_i,x_j) 
+ \sum_{u \in N(i)\setminus j} J^{(t)}_{u \rightarrow i}(y_i)\right)
+ \kappa^{(t+1)}_{i\ra j},
\end{equation}
where running estimates of the optimal solution are generated according to
\begin{equation}\label{eq:quadratic-minimizer-update}
\tilde{x}_i^{(t+1)} = \argmin_{y_i} \left(\tilde{f}^{(t+1)}_i(y_i) + \sum_{u \in N(i)} J^{(t+1)}_{u \rightarrow i}(y_i)\right).
\end{equation}
Note that the message update equation \eqref{eq:quadratic-update}
takes the form of a Ricatti equation for a scalar system, which can be
carried out efficiently.  Further, each optimization problem
\eqref{eq:quadratic-minimizer-update} is a scalar unconstrained convex
quadratic program.

A second approach makes use of a piecewise-linear approximation to
each message.  Let us assume knowledge that the optimal solution $x^*$
is in a closed bounded set $[-B,B]^n$.  Let $\Sscr =
\{\hat{x}_1,\ldots, \hat{x}_m\} \subset [-B,B]$, with $-B = \hat{x}_1
< \cdots < \hat{x}_m = B$, be a set of points where the linear pieces
begin and end.  Our approach applies the min-sum update equation to
compute values at these points.  Then, an approximation to the min-sum
message is constructed via linear interpolation between consecutive
points or extrapolation beyond the end points.  In particular, the
algorithm takes the form
\begin{equation}\label{eq:piecewise-update}
\begin{split}
J_{i \rightarrow j}^{(t+1)}(x_j) & = 
\min_{y_i \in [-B,B]} \left(f_i(y_i) + f_{ij}(y_i,x_j) 
+ \sum_{u \in N(i)\setminus j} J^{(t)}_{u \rightarrow i}(y_i)\right)
\\
& \quad\quad
+ \kappa^{(t+1)}_{i\ra j},
\end{split}
\end{equation}
for $x_j \in \Sscr$, where
\begin{equation}\label{eq:piecewise-update2}
J_{u \rightarrow i}^{(t)}(x_i) = \max_{1 \leq k \leq m-1} \frac{
(\hat{x}_{k+1} - x_i) J_{u \rightarrow i}^{(t)}(\hat{x}_{k+1}) + 
(x_i - \hat{x}_k) J_{u \rightarrow i}^{(t)}(\hat{x}_k)}
{\hat{x}_{k+1} - \hat{x}_k},
\end{equation}
for all $x_i \in \Re$.  As opposed to the case of quadratic
approximations, where each message is parameterized by two numerical
values, the number of parameters for each piecewise linear message
grows with $m$.  Hence, we anticipate that for fine-grain
approximations, our second approach is likely to require greater
computational resources.  On the other hand, piecewise linear
approximations may extend more effectively to non-convex problems,
since non-convex messages are unlikely to be well-approximated by
convex quadratic functions.

\section{Open Issues}\label{se:open}

There are many open questions in the theory of message passing
algorithms.  They fuel a growing research community that cuts across
communications, artificial intelligence, statistical physics,
theoretical computer science, and operations research.  This paper has
focused on application of the min-sum message passing algorithm to
convex programs, and even in this context a number of interesting
issues remain unresolved.

Our proof technique establishes convergence under total asynchronism
assuming a scaled diagonal dominance condition.  With such a flexible
model of asynchronous computation, convergence results for gradient
descent and coordinate descent also require similar diagonal dominance
assumptions.  On the other hand, for the {\it partially asynchronous}
setting, where communication delays and times between successive
updates are bounded, such assumptions are no longer required to
guarantee convergence of these two algorithms.  It would be
interesting to see whether convergence of the min-sum algorithm under
partial asynchronism can be established in the absence of scaled
diagonal dominance.

Another direction will be to assess practical value of the min-sum
algorithm for convex optimization problems.  This calls for
theoretical or empirical analysis of convergence and convergence times
for implementable variants as those proposed in the previous section.
Some convergence time results for a special case reported in
\cite{Moallemi06a} may provide a starting point.  Our expectation is
that for most relevant centralized optimization problems, the min-sum
algorithm will be more efficient than gradient descent or coordinate
descent but fall short of Newton's method.  On the other hand,
Newton's method does not decentralize gracefully, so in applications
that call for decentralized solution, the min-sum algorithm may prove
to be useful.

Finally, it would be interesting to explore whether ideas from this
paper can be helpful in analyzing behavior of the min-sum algorithm
for non-convex programs.  It is encouraging that convex optimization
theory has more broadly proved to be useful in designing and analyzing
approximation methods for non-convex programs.

\section*{Acknowledgments}

The first author was supported by a Benchmark Stanford Graduate
Fellowship.  This research was supported in part by the National
Science Foundation through grant IIS-0428868.

\bibliography{bp}

\end{document}